\documentclass[11pt]{article}

\usepackage{mathtools,amsmath,amssymb,color,amsthm}
\usepackage{mathrsfs}
\usepackage{authblk}
\usepackage{url}
\usepackage{hyperref}

\newtheorem{theorem}{Theorem}
\newtheorem{lemma}{Lemma}

\begin{document}

\title{An Elementary Proof of Sendov's Conjecture
via the Polar--Origin Inequalities}

\author{Jorge Alberto Borrego-Morell\\
\small Universidade Federal do Rio de Janeiro (UFRJ)\\
\small Campus Duque de Caxias Professor Geraldo Cidade\\
\small Rodovia Washington Luiz, 19593, km 104.5\\
\small Santa Cruz da Serra, Duque de Caxias, RJ 25240-005, Brazil\\
\small Corresponding author:
\href{mailto:jborrego@caxias.ufrj.br}{\texttt{jborrego@caxias.ufrj.br}}\\
\small ORCID: \href{https://orcid.org/0000-0001-8871-2739}
{0000-0001-8871-2739}}

\date{}
\maketitle

%

\begin{abstract}
In 1958, Blagovest Sendov conjectured that for every complex polynomial
$p$ of degree $n\geq2$ whose zeros all lie in the closed unit disk,
each zero of $p$ has at least one critical point within distance $1$.
We give a short elementary proof of Sendov's conjecture, establishing
the strict inequality at every interior zero. The principal ingredient
is a uniform analytical incompatibility estimate for the polar and
origin inequalities underlying the Mazur--Tao approach. By combining
exact binomial moment estimates with midpoint interpolation of
normalized excess moments, we eliminate the need for numerical
certification and obtain an argument valid in arbitrary degree.
The proof consists of a single real-variable lemma and its
application to polynomials. Combined with Rubinstein's classical
boundary theorem, our result also yields the Phelps--Rodriguez
equality characterization.
\end{abstract}



\section{Introduction}\label{In}

The relationship between the zeros of a polynomial and the zeros of
its derivatives is a central subject in the geometry of polynomials.
The classical Gauss--Lucas theorem describes the localization of the
critical points in the convex hull of the zeros, but it does not, by
itself, give a sufficiently precise estimate for the distance from an
individual zero to the set of critical points. The distinction between
these two questions underlies a substantial part of the analytic
theory of polynomials; see the monographs
\cite{Marden1966,RahmanSchmeisser2002,SheilSmall2002}.

Write $\mathbb D=\{z\in\mathbb C:|z|<1\}$ and, for a polynomial
$P$ of degree $n\geq2$ and a zero $a$ of $P$, set
\[
    \rho_a(P)=\min_{P'(w)=0}|w-a|.
\]
In 1958, Blagovest Sendov posed the conjecture that
\begin{equation}
\label{eq:intro-sendov}
    \rho_a(P)\leq1
\end{equation}
whenever all the zeros of $P$ belong to $\overline{\mathbb D}$.
The problem was subsequently recorded, under the name of Ilieff,
as Problem~4.5 in Hayman's collection of research problems
\cite{Hayman1967}; for its history and development, see
\cite{Marden1983,Sendov1995,Sendov2002}. The constant $1$ in
\eqref{eq:intro-sendov} is optimal for every degree, as is seen from
the polynomials $P(z)=c(z^n-\omega)$ with $c\ne0$ and $|\omega|=1$.
The strengthening proposed by Phelps and Rodriguez
\cite{PhelpsRodriguez1972} asserts that these polynomials are precisely
the possible equality configurations. In particular, it predicts the
strict inequality $\rho_a(P)<1$ at every interior zero $a$.

\medskip
\noindent\emph{Historical background.}
Progress on Sendov's conjecture occurred along several complementary
lines. The quadratic case is immediate, and the early low-degree
results of Brannan and Rubinstein \cite{Brannan1968,Rubinstein1968}
were followed by the proof for degree five by Meir and Sharma
\cite{MeirSharma1969}. Brown established the degree-six case
\cite{Brown1991}, while Borcea's investigations established the result
for polynomials with at most six and subsequently at most seven
distinct zeros \cite{Borcea1996Six,Borcea1996Seven}. The case of
degree seven was also treated by Brown \cite{Brown1997Seven}.
These developments culminated in the work of Brown and Xiang
\cite{BrownXiang1999}, which established the conjecture for degrees
$n\leq8$ and, more generally, for polynomials of arbitrary degree
having at most eight distinct zeros. For many years this remained the
established general low-degree range.

Another line of investigation concerned the geometry of the zero
configuration and the position of the distinguished zero. Rubinstein
\cite{Rubinstein1968} proved \eqref{eq:intro-sendov} when $|a|=1$,
including the strict inequality outside the exceptional family.
Further boundary-localization results were obtained by Goodman,
Rahman and Ratti \cite{GoodmanRahmanRatti1969} and by Joyal
\cite{Joyal1969}. Schmeisser \cite{Schmeisser1969} treated the case
in which the vertices of the convex hull of the zeros lie on the unit
circle. The general localization estimates of Bojanov, Rahman and
Szynal \cite{BojanovRahmanSzynal1985}, including the degree-dependent
bound $\rho_a(P)\leq2^{1/n}$, provided useful quantitative
approximations to the desired radius. Extremal configurations and
the associated optimization problems were investigated, in particular,
by Brown and Miller \cite{Brown1988,Miller1990}.

The position of $a$ relative to the unit circle proved particularly
important. Independently, V\^aj\^aitu and Zaharescu
\cite{VajaituZaharescu1993} and Miller \cite{Miller1993} established
the conjecture for a distinguished zero in a sufficiently thin
annulus adjacent to the unit circle, for each fixed degree. Subsequent
work refined the quantitative behavior of the Sendov radius near the
boundary \cite{Miller2005,Chijiwa2010,Chijiwa2011}. In particular,
Kasmalkar \cite{Kasmalkar2014} improved the quantitative near-boundary
range from an exponential to a polynomial scale in the degree.
These results, together with the elementary case $a=0$ and various
other special configurations, are part of the broader extremal theory
surveyed in \cite{Sendov2002,Schmeisser2002,Bojanov2011}.

A decisive development in the opposite direction, namely for large
degrees, was the work of D\'egot \cite{Degot2014}. For every fixed
$a\in(0,1)$, D\'egot proved the existence of a degree threshold
depending on $a$ above which Sendov's conjecture holds. Chalebgwa
\cite{Chalebgwa2020} made this dependence explicit and obtained
uniform thresholds when $a$ ranges over compact subintervals of
$(0,1)$. The dependence on the distinguished zero was removed by
Tao \cite{Tao2022}, who established the conjecture for all sufficiently
large degrees by an asymptotic compactness argument involving
logarithmic potentials, balayage and the argument principle. His
theorem did not supply an effective threshold separating the
large-degree regime from the remaining degrees. A subsequent preprint
of Zhang \cite{Zhang2026Effective} made the threshold explicit.

Besides these geometric and asymptotic approaches, the subject has
important connections with spectral theory, polynomial differentiators,
majorization and moment inequalities. We refer, for example, to
\cite{Pereira2003,Malamud2005,CheungNg2006,KhavinsonEtAl2011},
where related questions are formulated using matrices, compressions
or variances of zero configurations. These different formulations
illustrate the value of isolating scalar inequalities that express
the restrictions linking the zeros of a polynomial to its critical
points.

\medskip
\noindent\emph{The general-degree proof and its subsequent refinements.}
In August 2026, Mazur \cite{Mazur2026} released an AI-assisted,
computer-certified proof manuscript for Sendov's conjecture covering every
degree $n\geq2$, together with a separate Lean~4 formalization of the
theorem. The proof reduces a hypothetical counterexample to scalar
inequalities, whose terminal verification includes exact finite
certificates and a uniform auxiliary integral estimate. In particular,
the finite part of the manuscript treats $5\leq n-1\leq499$ through
16,862 certified interval records.

In his subsequent exposition, Tao \cite{Tao2026Digestion} identified
and organized the principal algebraic identities underlying Mazur's
argument, deriving what he calls the polar and origin inequalities.
The exposition also establishes the strict interior assertion and,
together with the boundary result, the full Phelps--Rodriguez
strengthening. A separate Lean development formalizes this argument
\cite{Tao2026Lean}. In the scalar analysis, Tao obtains the bound
$\alpha\leq17$ for a normalized parameter and then combines
large-degree analytic estimates with certified inequalities for the
remaining parameter range. The associated formalization records the
use of Bernstein positivity certificates. S.~An \cite{An2026} subsequently
reduced the computational part of this scalar argument by using
monotonicity, convexity and log-concavity, leaving a fixed collection
of explicit one-variable inequalities to be certified.

In a recent preprint,  Zhang
\cite{Zhang2026Companion} presents an analytical completion of the
Mazur--Tao argument based on the $D$-companion matrix. In place of
finite computational verification, his proof uses different integral
estimates in the two regimes $0<\alpha\leq3$ and $\alpha\geq3$,
where $\alpha=(n-1)(1-|a|^2)/2$.

Further developments have pursued genuinely stronger statements.
Tang and Zhang \cite{TangZhang2025} formulated a family of reciprocal
moment inequalities, and Zhang \cite{Zhang2026Quadratic} established
the quadratic inequality
\[
    \sum_{P'(w)=0}\frac{1}{|a-w|^2}\geq n-1,
\]
where critical points are counted with multiplicity and a zero
denominator is interpreted as giving $+\infty$. Recent preprints of
Bouthat and Fallat and of Zhang address the stronger quadratic
variance conjecture of Borcea
\cite{BouthatFallat2026,Zhang2026Variance,KhavinsonEtAl2011}.
A reciprocal-matrix treatment of Sendov's theorem and the quadratic
reciprocal inequality has also been developed by Tomov
\cite{Tomov2026}, while further conjectures concerning higher
derivatives have been treated by Ismoilov, Li, Rakhimov and Tan
\cite{IsmoilovEtAl2026}. Against this background, the present paper
focuses on the original scalar incompatibility problem in the
Mazur--Tao approach, giving a uniform analytical resolution directly
at the level of its two raw inequalities.

\medskip

\noindent\emph{The present contribution.}
We establish the following result
(Theorem~\ref{thm:strict-sendov-interior}):
if $P$ is a polynomial of degree $n\geq 2$ with all its zeros in
$\overline{\mathbb{D}}$ and $a$ is any zero of $P$ with $|a|<1$, then
\begin{equation}
\label{eq:intro-strict-sendov}
    \rho_a(P)<1.
\end{equation}
No additional restrictions are imposed on the location or
multiplicity of the other zeros. The case $a=0$ follows from
the product identity for $P'(0)$, while a rotation reduces
the remaining cases to $0<a<1$. Our proof is uniform in the
degree and requires no numerical verification.

The proof consists of an analytical incompatibility lemma
and its algebraic application to a polynomial. Following
\cite{Mazur2026,Tao2026Digestion}, assume $\rho_a(P)\geq 1$
and introduce the reciprocal critical-point coordinates
\[
    q_j=\frac{1}{a-w_j},\qquad
    x+iy=\frac{1}{m}\sum_{j=1}^{m}q_j,\qquad
    |q_j|\leq 1,
\]
where $w_1,\ldots,w_m$, with $m=n-1$, are the critical points
counted with multiplicity. Evaluations at the distinguished
zero, its reflected point $1/a$, and the origin, together
with the centroid identity, provide the algebraic reduction.
The historical antecedents of these identities include
\cite{Popoviciu1948,MeirSharma1969,Degot2014,
NazirMirWani2021,MirNazirWani2023};
see also \cite[Remark~8]{Tao2026Digestion}.

The algebraic reduction leads precisely to the two raw
inequalities of
\cite[Propositions~10(i) and~11(i)]{Tao2026Digestion}.
Writing
\[
    A=1-a^2,\qquad
    \alpha=\frac{mA}{2},\qquad
    \beta(t)=1-2axt+a^2t^2,
\]
and defining
\[
    G_m=a^2m(m+1)\int_0^1
    t\beta(t)^{(m-1)/2}\,dt,
\]
we obtain, for $m\geq 5$, the following necessary conditions
for a counterexample:
\begin{equation}
\label{eq:intro-raw-polar}
    1\leq\int_0^1
    (a^2+2aAtx+A^2t^2)^{m/2}\,dt
\end{equation}
and
\begin{equation}
\label{eq:intro-raw-origin}
    mA+ax\leq\frac{1-x^2}{2m}+G_m.
\end{equation}
The principal analytical contribution of this paper is to
establish their incompatibility directly in
Lemma~\ref{lem:analytical-polar-origin-incompatibility},
without requiring the supplementary bound $\alpha\leq 17$
or a finite covering of the scalar parameter domain.

The central simplification is obtained by setting
\[
    B=\beta(1),\qquad D=1-B.
\]
For $m\geq 5$, the polar inequality yields $0<B<1$ and
\begin{equation}
\label{eq:intro-polar-bound}
    mAD>6B,
\end{equation}
equivalently, $B<\alpha/(3+\alpha)$, the first endpoint
bound in Tao's polar estimate. In contrast to the separate
treatment of the two regimes of $\alpha$ in
\cite[Sections~3--4]{Zhang2026Companion}, we retain the raw
origin inequality \eqref{eq:intro-raw-origin} and establish
the uniform estimate
\begin{equation}
\label{eq:intro-uniform-origin}
    D(G_m-1)\leq\frac{24}{5}B
\end{equation}
for every $m\geq 5$ under the single polar condition
\eqref{eq:intro-polar-bound}.
This degree-independent estimate, valid throughout the
admissible parameter range, constitutes the essential
new analytical step.

To prove \eqref{eq:intro-uniform-origin}, the substitution
$v=at$ reduces $G_m$ to an integral involving powers of
$1-x^2+(x-v)^2$. Setting $r=x^2$ and $\delta=1-r$,
we control the contribution up to the minimum of this
quadratic by
\[
    L_m=m(m+1)\int_0^x
    v\bigl(\delta+(x-v)^2\bigr)^{(m-1)/2}\,dv.
\]
For odd $m\geq 5$, the binomial theorem and elementary
beta integrals yield $rL_m\leq 1+2\delta$.
Even degrees are treated by interpolating the excess moment
\[
    E_{\lambda}=\int_0^1(1-y)
    \bigl((\delta+ry^2)^{\lambda}-y^{2\lambda}\bigr)\,dy
\]
between neighboring integer exponents, rather than
interpolating the original integral. This preserves the
exact normalization at $x=1$, with an interpolation loss
controlled by
\[
    C_m=\sqrt{\frac{m^2-1}{m^2-4}}.
\]
When $a>x$, the polar constraint forces $mD>6$, allowing
the remaining contribution to be controlled by elementary
endpoint estimates.

Finally, \eqref{eq:intro-raw-origin} and
\eqref{eq:intro-uniform-origin} imply
\[
    \left(m-\frac{1}{2}\right)AD\leq\frac{27}{5}B,
\]
in contradiction with \eqref{eq:intro-polar-bound}.
The remaining degrees $m\leq 4$ are excluded by a single
elementary quartic integral identity.

The resulting argument simplifies the Mazur--Tao reduction
by establishing analytical incompatibility in arbitrary
degree using the binomial theorem, the Cauchy--Schwarz
inequality and elementary calculus, without computational
certification or a decomposition according to the size
of $\alpha$. Together with Rubinstein's classical boundary
theorem \cite{Rubinstein1968}, this yields a short elementary
proof of Sendov's conjecture and its Phelps--Rodriguez
equality characterization.

The analytical incompatibility lemma and the strict interior
theorem have also been formalized in Lean~4, verifying the
analytical estimates and their algebraic application
without resorting to numerical certificates.

\section{Proof of the results}

We first establish the scalar incompatibility on which the proof
rests. The following lemma collects all the real-variable estimates
needed in the argument and is formulated independently of the
polynomial setting.

\begin{lemma}[Analytical incompatibility of the polar and origin inequalities]
\label{lem:analytical-polar-origin-incompatibility}
Let $n\geq 2$ be an integer, $m=n-1$, $0<a<1$, and $x\in[-1,1]$.
Define
\[
    \alpha=\frac{m}{2}(1-a^2),
    \qquad
    \beta(t)=1-2axt+a^2t^2,\qquad 0\leq t\leq 1.
\]
Then the following two inequalities cannot hold simultaneously:
\begin{equation}
\label{eq:api-polar}
    1\leq\int_0^1
    \left(a^2+2at(1-a^2)x+t^2(1-a^2)^2\right)^{m/2}\,dt,
\end{equation}
and
\begin{equation}
\label{eq:api-origin}
    2\alpha+ax\leq\frac{1-x^2}{2m}
    +a^2nm\int_0^1t\,\beta(t)^{(m-1)/2}\,dt.
\end{equation}
\end{lemma}

\begin{proof}
Suppose that both inequalities hold. Put
\[
    A=1-a^2,\qquad B=\beta(1),\qquad D=1-B,
    \qquad Q(t)=a^2+2aAxt+A^2t^2.
\]
Thus $2\alpha=mA$. Define
\[
    G_m=a^2m(m+1)\int_0^1t\,\beta(t)^{(m-1)/2}\,dt.
\]
We shall prove the two estimates
\begin{equation}
\label{eq:api-combined-key}
    mAD>6B,
    \qquad
    D(G_m-1)\leq\frac{24}{5}B,
\end{equation}
from which the contradiction follows.

\smallskip
\noindent\emph{Polar estimate.}
The bases of all the real powers below are nonnegative: indeed,
for $0\leq t\leq1$,
\[
    Q(t)=(a-At)^2+2aAt(1+x)\geq0,\qquad
    Q(t)\leq(a+At)^2,
    \qquad \beta(t)=1-x^2+(x-at)^2\geq0.
\]
For $1\leq m\leq4$ and $v\geq0$, the weighted
arithmetic--geometric mean inequality gives
\[
    v^m\leq1-\frac m4+\frac m4v^4.
\]
Equivalently, multiplying the difference by $4$ gives, for
$m=1,2,3,4$, respectively,
\[
    (v-1)^2(v^2+2v+3),\quad 2(v^2-1)^2,\quad
    (v-1)^2(3v^2+2v+1),\quad 0.
\]
Direct integration gives
\[
\begin{aligned}
    1-\int_0^1(a+At)^4\,dt
    &=\frac{(1-a)^3}{5}
      (4+7a+3a^2-3a^3-2a^4+a^5)>0.
\end{aligned}
\]
The polynomial in parentheses is positive, since it equals
\[
    4+5a+3a^2(1-a)+2a(1-a)(1+a+a^2)+a^5.
\]
Consequently, if $m\leq4$,
\[
    1\leq\int_0^1 Q(t)^{m/2}\,dt
    \leq1-\frac m4+\frac m4\int_0^1(a+At)^4\,dt<1,
\]
a contradiction. Thus $m\geq5$. Since $B\geq(1-a)^2>0$ and
\[
    Q(t)=1-A+A(1+D)t-A^2t(1-t),
\]
inequality \eqref{eq:api-polar} forces $D>0$; otherwise the affine
majorant joining $Q(0)=a^2<1$ and $Q(1)=1+AD\leq1$ is strictly
smaller than $1$ for $t<1$. Hence $0<B,D<1$ and $x>a/2>0$.

For $0<t<1$, the affine majorant exceeds $Q(t)$ by
$A^2t(1-t)>0$. Since $m/2>0$, applying $1+u\leq e^u$ to
the majorant and integrating gives
\[
    1\leq\int_0^1 Q(t)^{m/2}\,dt
    <\int_0^1e^{\alpha[(1+D)t-1]}\,dt
    =e^{-\eta}\frac{\sinh h}{h},
\]
where $h=\alpha(1+D)/2$ and $\eta=\alpha B/2$.
We use the elementary inequality
\begin{equation}
\label{eq:api-combined-hyperbolic}
    \log\frac{\sinh u}{u}\leq\sqrt{u^2+9}-3,\qquad u>0.
\end{equation}
For a nonsingular proof, set $R=\sqrt{u^2+9}$ and
$Y(u)=u e^{R-3}$ for $u\geq0$. We have $Y(0)=0$,
$Y'(0)=1$, $R\geq3$, and direct differentiation gives
\[
    Y''(u)-Y(u)=e^{R-3}\,
        \frac{u(R-3)(2R-3)}{R^3}\geq0.
\]
Thus $K(u)=e^u(Y'(u)-Y(u))$ satisfies $K'(u)\geq0$ and
$K(0)=1$, so $Y'(u)-Y(u)\geq e^{-u}$. It follows that
\[
    H(u)=e^{-u}Y(u)+\frac12e^{-2u}
\]
has nonnegative derivative and $H(0)=1/2$. Therefore
$e^{-u}Y(u)\geq(1-e^{-2u})/2$, or $Y(u)\geq\sinh u$.
For $u>0$, division by $u$ and taking logarithms give
\eqref{eq:api-combined-hyperbolic}. This comparison uses only
regular functions at $u=0$.

It follows that $\eta<\sqrt{h^2+9}-3$, or $\eta^2+6\eta<h^2$.
Since $h^2-\eta^2=\alpha^2D$, we deduce the first estimate in
\eqref{eq:api-combined-key}:
\begin{equation}
\label{eq:api-combined-polar}
    mAD>6B.
\end{equation}

\smallskip
\noindent\emph{Origin estimate: the integral up to its minimum.}
Put $r=x^2$, $\delta=1-r$, and $p=(m-1)/2$. Note that
\begin{equation}
\label{eq:api-combined-geometry}
    D=r-(a-x)^2,\qquad B=\delta+(a-x)^2,
    \qquad 0<D\leq r,\quad 0\leq\delta\leq B.
\end{equation}
With the substitution $v=at$, we have
\[
    G_m=m(m+1)\int_0^a v\bigl(\delta+(x-v)^2\bigr)^p\,dv.
\]
For each integer $j\geq5$, define
\[
    L_j=j(j+1)\int_0^x
      v\bigl(\delta+(x-v)^2\bigr)^{(j-1)/2}\,dv.
\]

First let $j=2k+1$ be odd, with $k\geq2$. The binomial theorem gives
\[
    rL_j=r^2(2k+1)(2k+2)
      \sum_{i=0}^k\binom{k}{i}
      \frac{r^i\delta^{k-i}}{(2i+1)(2i+2)}.
\]
If $W_\ell=\binom{k+2}{\ell}r^\ell\delta^{k+2-\ell}$, this identity
becomes
\[
    rL_j=\sum_{i=0}^k
      \left(1+\frac{3(k-i)}{(2i+1)(k+2)}\right)W_{i+2}.
\]
The binomial theorem gives
$\sum_{\ell=0}^{k+2}W_\ell=1$ and, using
$(k+2-\ell)\binom{k+2}{\ell}
=(k+2)\binom{k+1}{\ell}$ for $0\leq\ell\leq k+1$,
$\sum_\ell(k+2-\ell)W_\ell=(k+2)\delta$.
Since $3/(2i+1)\leq1$ for $i\geq1$, these identities imply
\[
    rL_j\leq1+\delta+\frac{2k}{k+2}W_2
       =1+\delta+k(k+1)r^2\delta^k.
\]
To bound the last term without optimization, put
$V_\ell=\binom{k+1}{\ell}r^\ell\delta^{k+1-\ell}$.
For $k\geq2$, the first three positive-index weights satisfy
\[
\begin{aligned}
    V_1+V_3-V_2
    &=(k+1)r\delta^{k-2}
       \left(\delta^2-\frac{k}{2}r\delta
                    +\frac{k(k-1)}6r^2\right)\\
    &=(k+1)r\delta^{k-2}
       \left[\left(\delta-\frac{k}{4}r\right)^2
                    +\frac{k(5k-8)}{48}r^2\right]\geq0.
\end{aligned}
\]
Since the $V_\ell$ are nonnegative and sum to one, it follows that
\[
    k(k+1)r^2\delta^{k-1}=2V_2\leq V_1+V_2+V_3\leq1.
\]
This argument also covers $\delta=0$ (with the usual convention
that a zeroth power is $1$).
Consequently,
\begin{equation}
\label{eq:api-combined-odd}
    rL_j\leq1+2\delta\qquad(j\geq5\text{ odd}).
\end{equation}

Now suppose that $m\geq6$ is even. For $j\geq5$ put
$p_j=(j-1)/2$ and define
\[
    q(y)=\delta+ry^2,\qquad
    E_\lambda=\int_0^1(1-y)\bigl(q(y)^\lambda-y^{2\lambda}\bigr)\,dy,
    \qquad
    \mathcal E_j=r^2j(j+1)E_{p_j}.
\]
Here the exponents used in $E_\lambda$ satisfy $\lambda\geq2$.
Substituting $y=1-v/x$ and using
$\int_0^1(1-y)y^{j-1}\,dy=1/[j(j+1)]$ gives
\[
    rL_j=r^2+\mathcal E_j.
\]
Because $q(y)=y^2+\delta(1-y^2)\geq y^2$ for $0\leq y\leq1$,
we may express the normalized excess in terms of the densities
\[
    d_j(u)=j(j+1)p_j u^{p_j-1},\qquad u\geq0,
\]
as
\[
    \mathcal E_j=r^2\int_0^1(1-y)
                \int_{y^2}^{q(y)} d_j(u)\,du\,dy.
\]
The densities are nonnegative, including at $u=0$.
Since $p_m=(p_{m-1}+p_{m+1})/2$, direct calculation yields
\[
    d_m(u)^2=C_m^2 d_{m-1}(u)d_{m+1}(u),
    \qquad
    C_m=\sqrt{\frac{m^2-1}{m^2-4}}.
\]
The arithmetic--geometric mean inequality hence implies
\[
    d_m(u)\leq\frac{C_m}{2}\bigl(d_{m-1}(u)+d_{m+1}(u)\bigr).
\]
Integrating against the nonnegative weight above, we obtain
\[
    \mathcal E_m\leq\frac{C_m}{2}
                  \bigl(\mathcal E_{m-1}+\mathcal E_{m+1}\bigr).
\]
For the neighboring odd indices, \eqref{eq:api-combined-odd}
gives $\mathcal E_{m-1},\mathcal E_{m+1}\leq
1+2\delta-r^2=4\delta-\delta^2$. Consequently,
\[
    rL_m\leq r^2+C_m(4\delta-\delta^2).
\]
Since $1\leq C_m\leq1+3/[2(m^2-4)]$, we obtain
\[
    r(L_m-1)\leq(4C_m-1)\delta
       \leq\left(3+\frac6{m^2-4}\right)\delta.
\]
Thus both parities give $r(L_m-1)\leq
(3+\varepsilon_m)\delta$, where $\varepsilon_m$ is as below.
To pass from $r$ to $D$, if $L_m\leq1$ the desired estimate
is immediate from $D>0$ and $B>0$. If $L_m>1$, multiply
$D\leq r$ by the positive factor $L_m-1$ and use $\delta\leq B$.
In either case,
\begin{equation}
\label{eq:api-combined-left}
    D(L_m-1)\leq(3+\varepsilon_m)B,
    \qquad
    \varepsilon_m=
      \begin{cases}
        0,&m\text{ odd},\\[2pt]
        6/(m^2-4),&m\text{ even}.
      \end{cases}
\end{equation}

\smallskip
\noindent\emph{Origin estimate: the remaining integral.}
If $a\leq x$, then $G_m\leq L_m$, so \eqref{eq:api-combined-left}
already gives $D(G_m-1)<24B/5$. Suppose now that $a>x$ and put
\[
    z=\frac{a-x}{x}\in(0,1),\qquad T=G_m-L_m.
\]
By \eqref{eq:api-combined-geometry}, $B=\delta+rz^2$ and
$D=r(1-z^2)$. Also
\[
    T=rm(m+1)\int_0^z(1+w)(\delta+rw^2)^p\,dw.
\]
The bound $\delta+rw^2\leq B$ gives
\begin{equation}
\label{eq:api-combined-tail1}
    \frac{DT}{B}\leq
    m(m+1)D^2B^{p-1}\frac{z+z^2/2}{1-z^2}.
\end{equation}
Alternatively, using $w^2\leq zw$ in the term without the factor $w$
and integrating the term with $w$ exactly, we get
\[
    T\leq m\left(1+\frac2z\right)(B^{p+1}-\delta^{p+1}),
\]
and hence
\begin{equation}
\label{eq:api-combined-tail2}
    \frac{DT}{B}\leq mDB^p\left(1+\frac2z\right).
\end{equation}
These estimates apply to both integer and half-integer $p$.

Since $a>x$, we have $B-A=2a(a-x)>0$, and
\eqref{eq:api-combined-polar} yields $\tau:=mD>6$.
As $D<1$ and $m$ is an integer, $m\geq7$.
For $0<B\leq1$,
\[
    -\log B\geq\frac{2(1-B)}{1+B}=\frac{2D}{2-D};
\]
indeed, the difference of the two sides as functions of $B$
vanishes at $1$ and has derivative
$-(1-B)^2/[B(1+B)^2]\leq0$. Since $p-1=(m-3)/2$ and
$mD\geq6$, we obtain
\[
    (p-1)(-\log B)\geq\frac{(m-3)D}{2-D}
                      \geq\frac{\tau}{2}.
\]
The last inequality follows by multiplying by $2-D>0$ and
using $D(mD-6)\geq0$. In particular,
$e^{\tau/2}B^{p-1}\leq1$.

We record rational bounds for the resulting decay. Set
$y=(\tau-6)/4\geq0$. Then $\tau/6\leq1+y\leq e^y$, so
\[
    e^{\tau/2}=e^3e^{2y}\geq20(\tau/6)^2,
\]
where $e^3>20$ follows already from the terms through degree $8$
in the exponential series at $3$. Combining with the preceding
bound and using $0<B<1$ and $\tau\geq6$ gives
\[
    \tau^2B^{p-1}\leq\frac95,\qquad
    \tau B^p=\frac{B}{\tau}\,\tau^2B^{p-1}\leq\frac3{10}.
\]
Hence if $z\leq1/2$, \eqref{eq:api-combined-tail1} and
$(z+z^2/2)/(1-z^2)\leq5/6$ give
\[
    \frac{DT}{B}\leq
    \frac56\left(1+\frac1m\right)\frac95
    =\frac32\left(1+\frac1m\right).
\]
If $z\geq1/2$, \eqref{eq:api-combined-tail2} and
$1+2/z\leq5$ instead give $DT/B\leq3/2$, which is no larger.
In view of \eqref{eq:api-combined-left}, in both cases we have
\[
    \frac{D(G_m-1)}{B}
    \leq3+\varepsilon_m+\frac32\left(1+\frac1m\right).
\]
For odd $m\geq7$, this upper bound is at most $3+12/7<24/5$;
for even $m\geq8$, it is at most
$3+1/10+27/16<24/5$. We have thus established, in every case,
\begin{equation}
\label{eq:api-combined-origin}
    D(G_m-1)\leq\frac{24}{5}B,
\end{equation}
which completes the proof of \eqref{eq:api-combined-key}.

\smallskip
\noindent\emph{Final contradiction.}
Since $ax=1-(A+B)/2$, the origin inequality \eqref{eq:api-origin}
implies
\[
    \left(m-\frac12\right)AD
    \leq D(G_m-1)+\frac D2\left(B+\frac{\delta}{m}\right).
\]
Using \eqref{eq:api-combined-origin}, $D<1$, $\delta\leq B$, and
$m\geq5$, we obtain
\[
    \left(m-\frac12\right)AD
    \leq\left(\frac{24}{5}+\frac12+\frac1{2m}\right)B
    \leq\frac{27}{5}B.
\]
On the other hand, \eqref{eq:api-combined-polar} gives
\[
    \left(m-\frac12\right)AD
    >6\left(1-\frac1{2m}\right)B\geq\frac{27}{5}B,
\]
a contradiction. Therefore \eqref{eq:api-polar} and
\eqref{eq:api-origin} cannot hold simultaneously.
\end{proof}

The analytical incompatibility having been established, we now
return to the polynomial setting. We show that any hypothetical
counterexample to the strict Sendov bound necessarily satisfies
both inequalities of
Lemma~\ref{lem:analytical-polar-origin-incompatibility},
thereby obtaining the following result.

\begin{theorem}[Strict Sendov bound for an interior zero]
\label{thm:strict-sendov-interior}
Let $P$ be a polynomial of degree $n\geq2$, all of whose zeros
belong to $\overline{\mathbb D}$. If $a\in\mathbb D$ is a zero
of $P$, then
\[
    \rho_a(P)<1.
\]
\end{theorem}

\begin{proof}

Put $m=n-1$. Since multiplication by a nonzero constant does
not affect the critical distances, we may assume that $P$ is
monic and write
\[
    P(z)=(z-a)\prod_{j=1}^{m}(z-z_j),
    \qquad |z_j|\leq1,\quad 1\leq j\leq m.
\]
The remaining zeros $z_1,\ldots,z_m$ and the critical points
$w_1,\ldots,w_m$ are enumerated with multiplicity. No distinctness
hypothesis is imposed. If $a$ is a multiple zero, $P'(a)=0$
already gives the conclusion. If $a=0$, comparison of the two
factorizations of $P'(0)$ gives
\[
    n\prod_{j=1}^m|w_j|=\prod_{j=1}^m|z_j|\leq1,
\]
so that $\rho_0(P)\leq n^{-1/m}<1$. If $a\ne0$, put
$\lambda=a/|a|$. The polynomial $P(\lambda z)$ has a zero
at $|a|$ and its derivative is $\lambda P'(\lambda z)$.
Thus rotation, followed if necessary by multiplication by a
nonzero constant to restore monicity, preserves the degree,
the disk condition on zeros, and all critical distances.
We may therefore assume that $0<a<1$.

Suppose, to the contrary, that $\rho_a(P)\geq1$. In particular,
$P'(a)\ne0$ and $z_j\ne a$ for all $j$. Define
\[
    q_j=\frac1{a-w_j},\qquad 0<|q_j|\leq1,\qquad
    u=x+iy=\frac1m\sum_{j=1}^m q_j.
\]
The triangle inequality gives $|u|\leq m^{-1}\sum_j|q_j|\leq1$,
so $x^2+y^2\leq1$. Put
\[
    A=1-a^2,\qquad \alpha=\frac{mA}{2},\qquad
    \beta(t)=1-2axt+a^2t^2\quad(0\leq t\leq1).
\]
We will derive the two incompatible inequalities in
Lemma~\ref{lem:analytical-polar-origin-incompatibility}.

\medskip
\noindent\emph{ The polar inequality.}
Since
\[
    P'(a)=n\prod_{j=1}^m(a-w_j)
      =\frac{n}{\prod_{j=1}^m q_j}
      =\prod_{j=1}^m(a-z_j),
\]
we have
\[
    P'(z)=P'(a)\prod_{j=1}^m\bigl(1+(z-a)q_j\bigr).
\]
Integration along the segment from $a$ to $z$ yields
\begin{equation}
\label{eq:sspo-integral-representation}
    P(z)=(z-a)P'(a)\int_0^1
      \prod_{j=1}^m\bigl(1+t(z-a)q_j\bigr)\,dt.
\end{equation}
Evaluating at $z=1/a$, comparing with the factorization of $P$,
and multiplying by $a^m$, we obtain
\begin{equation}
\label{eq:sspo-polar-identity}
    \prod_{j=1}^m\frac{1-az_j}{a-z_j}
      =\int_0^1\prod_{j=1}^m(a+tAq_j)\,dt.
\end{equation}
For every $j$, the condition $|z_j|\leq1$ gives
\[
    |1-az_j|^2-|a-z_j|^2=A(1-|z_j|^2)\geq0.
\]
Consequently, by the triangle inequality, the arithmetic--geometric
mean inequality, and $|q_j|\leq1$,
\[
\begin{aligned}
    1&\leq\int_0^1\prod_{j=1}^m|a+tAq_j|\,dt\\
     &\leq\int_0^1
        \left(\frac1m\sum_{j=1}^m|a+tAq_j|^2\right)^{m/2}\,dt\\
     &\leq\int_0^1(a^2+2aAtx+A^2t^2)^{m/2}\,dt.
\end{aligned}
\]
In particular,
\begin{equation}
\label{eq:sspo-polar-inequality}
    1\leq\int_0^1
       (a^2+2aAtx+A^2t^2)^{m/2}\,dt.
\end{equation}
This is precisely \eqref{eq:api-polar}. Notice that only the closed
unit disk condition on the zeros has been used.

The small-degree argument at the start of the proof of
Lemma~\ref{lem:analytical-polar-origin-incompatibility} uses only
the polar inequality and excludes $1\leq m\leq4$.
We may henceforth assume $m\geq5$.

\medskip
\noindent\emph{Origin identities and the centroid identity.}
Let
\[
    F(t)=\prod_{j=1}^m(1-atq_j),\qquad
    I=\int_0^1F(t)\,dt,\qquad
    J=\prod_{j=1}^m z_jq_j,\qquad r=|J|\leq1,
    \qquad
    J_j=\left(\prod_{k\ne j}z_k\right)
        \left(\prod_{k=1}^m q_k\right).
\]
In particular $J=z_jJ_j$, including when $z_j=0$.
Evaluating \eqref{eq:sspo-integral-representation} at the origin gives
\begin{equation}
\label{eq:sspo-first-origin-identity}
    I=\frac{(-1)^mJ}{n}.
\end{equation}
Similarly, evaluating $P'(0)$ using the factorizations of $P$ and $P'$
gives
\begin{equation}
\label{eq:sspo-second-origin-identity}
    F(1)=\frac{(-1)^m}{n}
       \left(J+a\sum_{j=1}^m J_j\right).
\end{equation}
This is a polynomial cofactor identity and does not involve
division by any $z_j$.

Comparing the coefficients of the two highest powers in $P$ and $P'$
gives the centroid identity
\[
    \frac{a+\sum_{j=1}^m z_j}{n}
     =\frac1m\sum_{j=1}^m w_j
     =a-\frac1m\sum_{j=1}^m\frac1{q_j}.
\]
Taking conjugates and rearranging,
\begin{equation}
\label{eq:sspo-centroid-identity}
    \sum_{j=1}^m\frac1{\overline{q_j}}
       +\frac mn\sum_{j=1}^m\overline{z_j}
       =\frac{am^2}{n}.
\end{equation}

We use the following elementary product defect estimate:
\begin{equation}
\label{eq:sspo-defect-inequality}
    \sum_{\ell=1}^N(1-s_\ell^2)\prod_{k\ne\ell}s_k
        \leq1-\prod_{\ell=1}^N s_\ell^2,
    \qquad 0\leq s_\ell\leq1.
\end{equation}
For completeness, the two-factor inequality
$s(1-t^2)+t(1-s^2)\leq1-s^2t^2$ follows because the
difference between the right- and left-hand sides is
$(1-s)(1-t)(1-st)\geq0$. The case $N=1$ is equality.
For the induction step put $S=\prod_{k=1}^{N-1}s_k$.
By the induction hypothesis the left-hand side of
\eqref{eq:sspo-defect-inequality} is bounded above by
$s_N(1-S^2)+(1-s_N^2)S$, which is at most $1-s_N^2S^2$.
Applying \eqref{eq:sspo-defect-inequality} to
$|z_1|,\ldots,|z_m|,|q_1|,\ldots,|q_m|$, using
\[
\begin{aligned}
    |J_j-J\overline{z_j}|&=(1-|z_j|^2)|J_j|,\\
    \left|Jq_j-\frac{J}{\overline{q_j}}\right|
        &=(1-|q_j|^2)
          \left(\prod_{k=1}^m|z_k|\right)
          \left(\prod_{k\ne j}|q_k|\right),
\end{aligned}
\]
yields
\[
    \sum_{j=1}^m|J_j-J\overline{z_j}|
       +\sum_{j=1}^m\left|Jq_j-\frac{J}{\overline{q_j}}\right|
       \leq1-r^2.
\]
Multiplying \eqref{eq:sspo-centroid-identity} by $J$ and using $m/n\leq1$,
we deduce that
\[
    \left|\frac mn\sum_{j=1}^m J_j
       +J\sum_{j=1}^m q_j-\frac{am^2}{n}J\right|
       \leq1-r^2.
\]
Since $\sum_{j=1}^m q_j=mu$, it follows that
\begin{equation}
\label{eq:sspo-J-sum}
    \sum_{j=1}^m J_j=(am-nu)J+E,
       \qquad |E|\leq\frac nm(1-r^2).
\end{equation}

\medskip
\noindent\emph{ The origin inequality.}
Combining \eqref{eq:sspo-first-origin-identity},
\eqref{eq:sspo-second-origin-identity}, and \eqref{eq:sspo-J-sum} gives
\[
    F(1)+amuI=\frac{(-1)^m}{n}
       \bigl[J(1+ma^2-au)+aE\bigr].
\]
Thus, setting $C=|1+ma^2-au|/n$, we obtain
\[
    |F(1)+amuI|\leq rC+\frac am(1-r^2).
\]
Since $m\geq5$ and $a>0$, the following square decomposition
gives the lower bound needed to absorb the product defect:
\[
\begin{aligned}
    m(1+ma^2-a)-2an
      &=(ma-2)^2+(m-4)+(m-2)a>0.
\end{aligned}
\]
Because $|u|\leq1$, we have $C\geq(1+ma^2-a)/n>2a/m$. Therefore
\[
    rC+\frac am(1-r^2)
       =C-(1-r)\left(C-\frac am(1+r)\right)\leq C.
\]
Moreover, $1+ma^2-ax\geq ma^2>0$ and $y^2\leq1-x^2$, whence
\[
\begin{aligned}
    nC&=\sqrt{(1+ma^2-ax)^2+a^2y^2}\\
       &\leq 1+ma^2-ax+\frac{a^2y^2}{2(1+ma^2-ax)}\\
       &\leq 1+ma^2-ax+\frac{1-x^2}{2m}.
\end{aligned}
\]
It follows that
\begin{equation}
\label{eq:sspo-origin-linear-bound}
    |F(1)+amuI|\leq\frac{1+ma^2-ax}{n}
                          +\frac{1-x^2}{2nm}.
\end{equation}

On the other hand, differentiation of $F$ gives the identity
\[
    F'(t)=-amuF(t)-a^2t\sum_{j=1}^m q_j^2
                               \prod_{k\ne j}(1-atq_k).
\]
Integrating and using $F(0)=1$ and $|q_j|\leq1$, we obtain
\[
    1\leq |F(1)+amuI|+a^2\int_0^1t
                  \sum_{j=1}^m\prod_{k\ne j}|1-atq_k|\,dt.
\]
The needed special case of Maclaurin's inequality can be obtained
directly. For nonnegative $v_1,\ldots,v_m$ of mean $M$, one has
$\sum_j\prod_{k\ne j}v_k\leq mM^{m-1}$.
The case $m=1$ is equality. For the induction step let $b$ be the
mean of the first $m-1$ entries and let $c=v_m$, so
$mM=(m-1)b+c$. The induction hypothesis and the arithmetic--geometric
mean inequality bound the omitted-product sum by
\[
    b^{m-1}+c(m-1)b^{m-2}
    =m\bigl((m-1)b^{m-2}M-(m-2)b^{m-1}\bigr)
    \leq mM^{m-1},
\]
where the final inequality is the tangent inequality for the
convex power $s\mapsto s^{m-1}$; for $m=2$ it is equality.
Applying this estimate and the Cauchy--Schwarz inequality,
\[
\begin{aligned}
    \sum_{j=1}^m\prod_{k\ne j}|1-atq_k|
      &\leq m\left(\frac1m\sum_{k=1}^m|1-atq_k|\right)^{m-1}\\
      &\leq m\left(\frac1m\sum_{k=1}^m|1-atq_k|^2\right)^{(m-1)/2}\\
      &\leq m\,\beta(t)^{(m-1)/2},
\end{aligned}
\]
where the last inequality uses $|q_k|\leq1$ and
$m^{-1}\sum_{k=1}^m\operatorname{Re}q_k=x$.
Combining these estimates with \eqref{eq:sspo-origin-linear-bound}
yields
\[
    1\leq\frac{1+ma^2-ax}{n}+\frac{1-x^2}{2nm}
         +a^2m\int_0^1t\,\beta(t)^{(m-1)/2}\,dt.
\]
Multiplying by $n$ and rearranging gives
\begin{equation}
\label{eq:sspo-origin-inequality}
    2\alpha+ax=m(1-a^2)+ax
       \leq\frac{1-x^2}{2m}
          +a^2nm\int_0^1t\,\beta(t)^{(m-1)/2}\,dt.
\end{equation}
This is precisely \eqref{eq:api-origin}.

\medskip
\noindent\emph{Conclusion.}
We have obtained both \eqref{eq:sspo-polar-inequality} and
\eqref{eq:sspo-origin-inequality}, with $n=m+1$, $0<a<1$, and
$x\in[-1,1]$. This contradicts
Lemma~\ref{lem:analytical-polar-origin-incompatibility}. Hence
$\rho_a(P)<1$.
\end{proof}

\medskip
 
Theorem~\ref{thm:strict-sendov-interior} establishes the strict
Sendov bound for every interior zero and every degree.
For a zero $a$ on the unit circle, Rubinstein's theorem
\cite[Theorem~1]{Rubinstein1968} establishes the same strict bound
except for polynomials of the form
\[
    P(z)=c(z^n-a^n),\qquad c\ne0,\quad |a|=1.
\]
In this exceptional case, all the critical points are at the origin
and hence $\rho_a(P)=1$.

Consequently, for every polynomial of degree $n\geq2$ with all
its zeros in $\overline{\mathbb D}$,
\[
    \rho_a(P)\leq1
\]
at every zero $a$, with equality precisely in the exceptional
configuration described above. This completes the elementary
proof of Sendov's conjecture and the Phelps--Rodriguez
equality characterization.

\section*{Formal verification}

Lemma~\ref{lem:analytical-polar-origin-incompatibility} and
Theorem~\ref{thm:strict-sendov-interior} have been formalized
and verified in Lean~4 using Mathlib. The formal development
follows the analytical proof and the subsequent algebraic
reduction presented in this paper. In particular, the final
theorem is stated for arbitrary complex polynomials and
arbitrary interior zeros.

The principal Lean declarations are
\begin{verbatim}
SendovProof.analytical_polar_origin_incompatibility
SendovProof.analytical_polar_origin_incompatibility_raw
SendovProof.strict_sendov_interior
\end{verbatim}
The first two are contained in \verb|PolarOrigin.lean|,
and the third in \verb|Sendov.lean|.

The development was checked using Lean~4.24.0 and Mathlib
at the revision
\begin{center}
\texttt{f897ebcf72cd16f89ab4577d0c826cd14afaafc7}.
\end{center}
For each of the three declarations, the transitive axiom
audit reports precisely
\begin{verbatim}
[propext, Classical.choice, Quot.sound]
\end{verbatim}
which are standard foundational axioms in Lean. In particular,
none of the verified results depends on \verb|sorryAx| or
any additional axioms. The formalization thus provides a
machine-checked verification of the analytical argument,
rather than supplying computational certificates needed
to complete the proof.

The formal verification concerns the strict Sendov inequality
at interior zeros. The boundary case is supplied by
Rubinstein's classical theorem
\cite[Theorem~1]{Rubinstein1968}, which, together with the
strict interior inequality, yields the full Phelps--Rodriguez
equality characterization. The boundary argument is not part
of the present formalization.

\paragraph{Code availability.}
The complete Lean~4 project, including the pinned toolchain,
dependency manifest, axiom-audit file, recorded audit results,
and instructions for reproducing the verification, is publicly
available at
\url{https://github.com/jbmorell/SendovProof}.

\section*{Acknowledgements}

The author acknowledges the use of ChatGPT (GPT-6.1 and GPT-6.2;
OpenAI) for mathematical exploration, the development and
examination of analytical arguments, manuscript preparation,
and assistance with the development and refinement of the
accompanying Lean~4 formalization.

The central idea of the proof and its conceptual organization,
whereby the main result is reduced to a single analytical lemma
and its application in the main theorem, originate from the author.

The author directed the AI-assisted research process, critically
reviewed and verified the generated output, and assumes full
responsibility for the mathematical results, the accompanying
formalization, and the final content of the manuscript.

\bibliography{sendov_references}

\end{document}